\documentclass[a4paper,reqno,12pt]{amsart}
\usepackage{amsmath,amsrefs,color}

\numberwithin{equation}{section}

\DeclareMathOperator{\Scal}{S}
\DeclareMathOperator{\Vol}{Vol}
\DeclareMathOperator{\OO}{O}
\DeclareMathOperator{\oo}{o}

\DeclareMathOperator{\Span}{span}
\newcommand{\R}{\mathbb{R}}

\newcommand{\<}{\left<}
\renewcommand{\>}{\right>}
\renewcommand{\[}{\left[}
\renewcommand{\]}{\right]}
\renewcommand{\(}{\left(}
\renewcommand{\)}{\right)}

\newtheorem{theorem}{Theorem}[section]

\newtheorem{proposition}[theorem]{Proposition}
\newtheorem{lemma}[theorem]{Lemma}

\begin{document}

\title[Positive clusters in dimensions four and five]{Positive clusters for smooth perturbations of a critical elliptic equation in dimensions four and five}

\author{Pierre-Damien Thizy}

\address{Pierre-Damien Thizy, Universit\'e de Cergy-Pontoise, CNRS, D\'epar\-tement de Math\'ematiques, F-95000 Cergy-Pontoise, France.}
\email{pierre-damien.thizy@u-cergy.fr}

\author{J\'er\^ome V\'etois}

\address{J\'er\^ome V\'etois, McGill University, Department of Mathematics and Statistics, 805 Sherbrooke Street West, Montreal, Quebec H3A 0B9, Canada.}
\email{jerome.vetois@mcgill.ca}

\date{March 22, 2016}

\begin{abstract}
We construct clustering positive solutions for a perturbed critical elliptic equation on a closed manifold of dimension $n=4,5$. Such a construction is already available in the literature in dimensions $n\ge 6$ (see for instance~\cites{ChenLin,DruHeb,MorPisVai,PisVai,RobVet3}) and not possible in dimension $3$ by~\cite{LiZhu}. This also provides new patterns for the Lin--Ni~\cite{LinNiCon} problem on closed manifolds and completes results by Br\'ezis and Li~\cite{BreLi} about this problem.
\end{abstract}

\maketitle

\section{Introduction and main result}\label{Sec1}

Let $(M^n,g)$ be a smooth closed Riemannian manifold of dimension $n\ge 3$, and $2^\star=\frac{2n}{n-2}$ be the critical Sobolev exponent for the embeddings of $H^1(M)$ into the Lebesgue spaces. Given smooth perturbations $(h_\varepsilon)_\varepsilon$ of a function $h_0$ in $M$, the asymptotic behavior of a sequence $(u_\varepsilon)_\varepsilon$ of smooth positive functions satisfying
\begin{equation}\label{Intro1}
\Delta_g u_\varepsilon+h_\varepsilon u_\varepsilon=u_\varepsilon^{2^\star-1}
\end{equation}
for all $\varepsilon>0$ has been intensively studied in the last decades. Here $\Delta_g=-\text{div}_g(\nabla\cdot)$ is the Laplace--Beltrami operator. If such a sequence $(u_\varepsilon)_\varepsilon$ is bounded in $H^1(M)$, then we know from Struwe~\cite{Struwe} that there exist $k\in\mathbb{N}$, $k$ sequences $(\mu_{1,\varepsilon})_\varepsilon,\dotsc,(\mu_{k,\varepsilon})_\varepsilon$ of positive numbers converging to~0, and $k$ sequences $(\xi_{1,\varepsilon})_\varepsilon,\dotsc,(\xi_{k,\varepsilon})_\varepsilon$ of points converging to $\xi_1,\dotsc,\xi_k$ in $M$ such that
\begin{equation}\label{Intro2}
u_\varepsilon=u_0+\sum_{i=1}^k \(\frac{\sqrt{n\(n-2\)}\mu_{i,\varepsilon}}{\mu_{i,\varepsilon}^2+d_g(\xi_{i,\varepsilon},\cdot)^2}\)^{\frac{n-2}{2}}+o(1),
\end{equation}
where $o(1)\to 0$ strongly and $u_\varepsilon\rightharpoonup u_0$ in $H^1(M)$ as $\varepsilon\to 0$. If the sequence $(u_\varepsilon)_\varepsilon$ is not uniformly bounded, then we say that $(u_\varepsilon)_\varepsilon$ blows up and in this case, it follows from classical elliptic estimates that $k$ is non-zero in \eqref{Intro2}. If $\xi_1=\dotsb=\xi_k=\xi_0$, then we say that $(u_\varepsilon)_\varepsilon$ blows up with $k$ peaks at the point $\xi_0$.

\smallskip
In the case of dimension~3, it was proved by Li and Zhu~\cite{LiZhu} (see Theorem~6.3 in Hebey~\cite{HebeyBook}) that $\xi_1,\dotsc,\xi_k$ are necessarily distinct in \eqref{Intro2}. By contrast, in the case of dimensions larger than or equal to~6, Druet and Hebey~\cite{DruHeb}, Robert and V\'etois~\cite{RobVet3}, and more recently, Morabito, Pistoia, and Vaira~\cite{MorPisVai} and Pistoia and Vaira~\cite{PisVai} have given examples of $(h_\varepsilon)_\varepsilon$ and $(u_\varepsilon)_\varepsilon$ for which $k\ge 2$ is arbitrary and the sequences $(\xi_{i,\varepsilon})_\varepsilon$, $i=1,\dotsc,k$, converge to the same point of $M$ in \eqref{Intro2} (see also Chen and Lin~\cite{ChenLin} where a similar result was obtained for the prescribed scalar curvature equation on the sphere in dimensions $n\ge7$). The main goal of this paper is to prove that such examples can actually be given starting from dimension $4$. We state our result as follows.

\begin{theorem}\label{Th}
Let $\(M,g\)$ be a closed manifold of dimension $n\in\left\{4,5\right\}$. Assume that the scalar curvature $\Scal_g$ of the manifold has a non-degenerate minimum point $\xi_0$ such that $\Scal_g\(\xi_0\)<0$. Then for any natural number $k>1$, there exists a family of positive solutions $\(u_{k,\varepsilon}\)_{\varepsilon>0}$ of the equations
\begin{equation}\label{Eq}
\Delta_gu_{k,\varepsilon}+\varepsilon u_{k,\varepsilon}=u_{k,\varepsilon}^{2^*-1}\quad\text{in }M
\end{equation}
such that $\(u_{k,\varepsilon}\)_{\varepsilon>0}$ blows up with $k$ peaks at the point $\xi_0$ as $\varepsilon\to0$.
\end{theorem}

According to the terminology of Schoen~\cite{Schoen}, the blow-up points of $(u_{\varepsilon,k})_\varepsilon$ that we construct in Theorem~\ref{Th} are non-isolated blow-up points. Isolation of blow-up points turned out to be a crucial step in the proofs of compactness for the Yamabe equation (see Druet~\cite{Dru}, Khuri, Marques, and Schoen~\cite{KhuMarSch}, Li and Zhang~\cites{LiZha1,LiZha2}, Li and Zhu~\cite{LiZhu}, \linebreak Marques~\cite{Mar}, Schoen~\cite{Sch2}, and Schoen and Zhang~\cite{SchZha}). Isolated blow-up points for the Yamabe equation were constructed by Brendle~\cite{Bre} and Brendle and Marques~\cite{BreMar} in high dimensions. However, as we explained above, with regards to sequences of solutions $(u_\varepsilon)_\varepsilon$ of more general perturbed critical elliptic equations like \eqref{Intro1}, the a priori blow-up analysis cannot rule out in general non-isolated blow-up points in dimensions $n\ge 4$ (see for instance~\cites{DruHeb,DruHebRob}). 

\smallskip
More specifically, looking for non-constant solutions of Equation \eqref{Eq} or ``patterns'' has some relevance in mathematical biology (see for instance~\cites{GieMei}). This is referred to in the literature as the Lin--Ni~\cite{LinNiCon} or Lin--Ni--Takagi~\cite{LinNiTak} problem. In the case of a closed manifold, Br\'ezis and Li~\cite{BreLi} proved that the only solution to \eqref{Eq} is the constant solution for $0<\varepsilon\ll1$. This result also holds true when the scalar curvature $\Scal_g$ is positive everywhere in dimensions $n\ge 4$ (see Druet~\cite{Dru} and Remark~6.1~$(i)$ in Hebey~\cite{HebeyBook}). The above Theorem~\ref{Th} proves that this result generically fails in dimensions $n=4,5$ when $\Scal_g$ is negative somewhere. Moreover the clustering solutions that we construct in Theorem~\ref{Th} give new type of patterns for the Lin--Ni problem. The first author proved in~\cite{Thi} that Theorem~\ref{Th} fails for all $k\ge1$ in dimensions $n\ge 6$ and also that in dimensions $n=4,5$, the non-degeneracy assumption in Theorem~\ref{Th} can be removed in case $k=1$ (see the one-peak version of Theorem~\ref{Th} in~\cite{Thi}).

\smallskip
There is an abundant literature about the original Lin--Ni problem, namely Equation \eqref{Eq} posed on a bounded domain of the Euclidean space with zero Neumann boundary condition. We mention of course the works of Lin and Ni~\cite{LinNiCon} and Lin, Ni, and Takagi~\cite{LinNiTak}, where after proving a subcritical analogue result in~\cite{LinNiTak}, it was conjectured in~\cite{LinNi} that this equation does not have any other solution than the constant solution for $0<\varepsilon\ll1$. Without any pretension to exhaustivity, we also mention Adimurthi and Yadava~\cites{LinNiAdiYad,LinNiAdiYad2} and Budd, Knapp, and Peletier~\cite{BudKnaPel} for a complete discussion of the radial case when the domain is a ball (conjecture false for $n=4,5,6$ and true otherwise) and Rey and Wei~\cite{ReyWei} and Wei, Xu, and Yang~\cite{WeiXuYang} who proved that the conjecture fails for all bounded domains of dimension $n=5$ and $n=4,6$ respectively. The solutions constructed in~\cite{ReyWei} and~\cite{WeiXuYang} have isolated blow-up points in the interior of the domain, one blow-up point in~\cite{WeiXuYang} and multiple blow-up points in~\cite{ReyWei}. Zhu~\cite{ZhuLinNi} proved that the Lin--Ni conjecture holds true in $3$--dimensional convex domains and Wang, Wei, and Yan~\cites{WanWeiYan,WanWeiYan2} proved that the conjecture fails for non-convex domains of dimension $n\ge 3$. Druet, Robert, and Wei~\cite{LinNiAMS} proved that the Lin--Ni conjecture is true in convex domains of dimension $n\not\in\{4,5,6\}$, assuming a bound on the energy of solutions. In the case where the parameter $\varepsilon$ does not approach zero, we mention for instance the works of del~Pino, Felmer, Rom\'an, and Wei~\cite{DelPMusRomWei} for $\varepsilon$ close to a fixed number, and Esposito~\cite{Esp}, Gui and Lin~\cite{GuiLin}, and Wei and Yan~\cite{WeiYan} for $\varepsilon$ converging to infinity, and we refer to these papers and the references therein for a more complete discussion. A vectorial version of the Lin--Ni conjecture has also been considered by Hebey~\cite{LinNi}.

\smallskip
The proof of Theorem~\ref{Th} relies on the Lyapunov--Schmidt method and uses the general formalism developped in Robert and V\'etois~\cite{RobVet2}. This allows to reduce the problem to finding critical points of an energy function on a finite dimensional space, here of dimension $k\(n+1\)+1$. In our case, we are dealing with a situation where the reduced energy function has a saddle point. To manage this type of situations, we prove a general critical point result in Appendix~A which allows to restrict the computations of $C^1$--estimates to a smaller number of variables. This generalizes an argument used by Chen, Wei, and Yan~\cite{ChenWeiYan} in the case of a function of two real variables. We believe this result may be useful in future works based on the Lyapunov--Schmidt method when dealing with a saddle point situation.

\smallskip
Another specificity of our constructions is the role played by the interaction between the peaks and the constant solutions. This can be seen by looking at the dependence on $\varepsilon$ of the parameter $z_{\varepsilon}$ in our approximated solutions, which are of the form
$$u_\varepsilon=z_\varepsilon+\sum_{i=1}^k \(\frac{\sqrt{n\(n-2\)}\mu_{i,\varepsilon}}{\mu_{i,\varepsilon}^2+d_g(\xi_{i,\varepsilon},\cdot)^2}\)^{\frac{n-2}{2}}+\phi_\varepsilon\,,$$
where $z_{\varepsilon}$ is a small positive parameter, $\mu_{i,\varepsilon}$ and $\xi_{i,\varepsilon}$ are as in \eqref{Intro2}, and $\phi_\varepsilon$ is a remainder term in $H^1\(M\)$ which is orthogonal to a finite dimensional subspace including the constant functions. While in dimension $n=5$, $z_{\varepsilon}$ behaves at first order like the constant solutions of \eqref{Eq}, namely $z_{\varepsilon}\sim\varepsilon^{3/4}$ as $\varepsilon\to0$, the situation becomes very different when the dimension $n$ jumps down to $4$ (see \eqref{z}). In this case, we find that $z_{\varepsilon}$ has exponential decay as $\varepsilon\to0$, which indicates that there is a much stronger interaction between this term and the peaks in dimension $n=4$, which as explained above, is the lowest possible dimension for the existence of positive clusters. This is also the reason why we obtain different expressions for the reduced energy functions in dimensions $n=4$ and $n=5$.

\smallskip
The paper is organized as follows. We introduce our ansatz of multipeak solutions and perform the main part of the proof of Theorem~\ref{Th} in Section~\ref{Sec2}. We perform the error estimates and $C^0$--energy estimates in Section~\ref{Sec3} and the $C^1$--energy estimates in Section~\ref{Sec4}. Finally we prove our general critical point result in Appendix~\ref{Appendix}.

\section{Proof of Theorem~\ref{Th}}\label{Sec2}

We fix $k>1$ and $\xi_0\in M$ as in the statement of Theorem~\ref{Th}. Since $M$ is compact, we may fix a positive real number $r_0$ such that $r_0$ is less than the injectivity radius at all points of the manifold $\(M,g\)$. For any real numbers $\varepsilon,K>0$, we consider the parameter set
\begin{multline*}
\mathcal{D}_{K,\varepsilon}:=\bigg\{\(\xi,\mu\)=\(\(\xi_1,\dotsc,\xi_k\),\(\mu_1,\dotsc,\mu_k\)\)\in B\(\xi_0,r_0\)^k\times\(0,\varepsilon\)^k:\\
\frac{\mu_i}{\mu_j}+\frac{\mu_j}{\mu_i}+\frac{d_g\(\xi_i,\xi_j\)^2}{\mu_i\mu_j}>K\quad\forall i\ne j\bigg\},
\end{multline*}
where $d_g\(\xi_i,\xi_j\)$ is the geodesic distance between $\xi_i$ and $\xi_j$, and $B\(\xi_0,r_0\)$ is the geodesic ball of center $\xi_0$ and radius $r_0$ in the manifold $\(M,g\)$. We let $\chi$  be a smooth cutoff function such that $0\le\chi\le1$ on $\[0,\infty\)$, $\chi=1$ on $\[0,r_0/2\]$, and $\chi=0$ on $\[r_0,\infty\)$. We consider the family of profiles
$$u_{z,\xi,\mu}\(x\):=z+\sum_{i=1}^kW_{\xi_i,\mu_i}\(x\),$$ 
for all $x\in M$ and $\(z,\xi,\mu\)\in\(0,\varepsilon\)\times\mathcal{D}_{K,\varepsilon}$, where 
$$W_{\xi_i,\mu_i}\(x\):=\chi\(d_g\(x,\xi_i\)\)\(\frac{\sqrt{n\(n-2\)}\mu_i}{\mu_i^2+d_g\(x,\xi_i\)^2}\)^{\frac{n-2}{2}}$$
for all $i\in\left\{1,\dotsc,k\right\}$.

\smallskip
For any real number $\varepsilon>0$, the energy functional of Equation \eqref{Eq} is defined as
$$J_\varepsilon\(u\)=\frac{1}{2}\int_M\(\left|\nabla u\right|_g^2+\varepsilon u^2\)dv_g-\frac{1}{2^*}\int_Mu_+^{2^*}dv_g$$
for all $u\in H^1\(M\)$, where $u_+:=\max\(u,0\)$. For any $\(z,\xi,\mu\)\in\(0,\varepsilon\)\times\mathcal{D}_{K,\varepsilon}$, we define our profile's error as
$$R_{\varepsilon,z,\xi,\mu}:=\left\|\(\Delta_g+\varepsilon\)u_{z,\xi,\mu}-u_{z,\xi,\mu}^{2^*-1}\right\|_{L^{\frac{2n}{n+2}}\(M\)}.$$

\smallskip
As a particular case of Theorem~1.1 of Robert and V\'etois~\cite{RobVet2}, we obtain the following result.

\begin{proposition}\label{Pr1}
There exist positive constants $\varepsilon_0$, $C_0$, and $K_0$, such that for any real number $\varepsilon\in\(0,\varepsilon_0\)$, there exists a mapping $\phi_\varepsilon\in C^1\(\(0,\varepsilon_0\)\times\mathcal{D}_{K_0,\varepsilon_0},H^1\(M\)\)$ such that for any $\(z,\xi,\mu\)\in\(0,\varepsilon_0\)\times\mathcal{D}_{K_0,\varepsilon_0}$, we have
\begin{align}
&\left|J_\varepsilon\(u_{z,\xi,\mu}+\phi_\varepsilon\(z,\xi,\mu\)\)-J_\varepsilon\(u_{z,\xi,\mu}\)\right|\le C_0\,R_{\varepsilon,z,\xi,\mu}^2\,,\label{Energy0}\\
&\left\|\phi_\varepsilon\(z,\xi,\mu\)\right\|_{H^1\(M\)}\le C_0\,R_{\varepsilon,z,\xi,\mu}\,,\label{Phi}
\end{align}
and 
\begin{multline}\label{Critical}
D_uJ_\varepsilon\(u_{z,\xi,\mu}+\phi_\varepsilon\(z,\xi,\mu\)\)=0\\
\Longleftrightarrow\(\partial_z\mathcal{J}_\varepsilon\(z,\xi,\mu\),D_\mu\mathcal{J}_\varepsilon\(z,\xi,\mu\),D_\xi\mathcal{J}_\varepsilon\(z,\xi,\mu\)\)=\(0,0,0\),
\end{multline}
where $\mathcal{J}_\varepsilon\(z,\xi,\mu\)=J_\varepsilon\(u_{z,\xi,\mu}+\phi_\varepsilon\(z,\xi,\mu\)\)$.
\end{proposition}

Now we need to specify the dependence of our parameters $\(z,\xi,\mu\)$ with respect to $\varepsilon$. For any $s\in\R$, $t=\(t_1,\dotsc,t_k\)\in\R^k$, and $\tau=\(\tau_1,\dotsc,\tau_k\)\in\(T_{\xi_0}M\)^k$, where $T_{\xi_0}M$ is the tangent space of $\(M,g\)$ at the point $\xi_0$, we define
\begin{equation}\label{z}
z_{\varepsilon,s}:=\left\{\begin{aligned}&\varepsilon^{-1}e^{-s/\varepsilon}&&\text{if }n=4\\&\varepsilon^{3/4}+s\,\varepsilon^{5/4}&&\text{if }n=5,\end{aligned}\right.
\end{equation}
$$\mu_{\varepsilon,s,t}=\(\mu_{\varepsilon,s,t_i}\)_{1\le i\le k}:=\(\mu_{\varepsilon,s}t_i\)_{1\le i\le k},\quad\mu_{\varepsilon,s}:=\left\{\begin{aligned}&e^{-s/\varepsilon}&&\text{if }n=4\\&\varepsilon^{3/2}&&\text{if }n=5,\end{aligned}\right.$$
and
$$\xi_{\varepsilon,\tau}=\(\xi_{\varepsilon,\tau_i}\)_{1\le i\le k}:=\(\exp_{\xi_0}\(\delta_\varepsilon\tau_i\)\)_{1\le i\le k},\quad\delta_\varepsilon:=\left\{\begin{aligned}&\varepsilon^{1/4}&&\text{if }n=4\\&\varepsilon^{3/10}&&\text{if }n=5.\end{aligned}\right.$$
In particular, we point out that for any $i,j\in\left\{1,\dotsc,k\right\}$, we have
\begin{equation}\label{dist}
d_g\(\xi_{\varepsilon,\tau_i},\xi_{\varepsilon,\tau_j}\)=\delta_\varepsilon\left|\tau_i-\tau_j\right|+\OO\(\delta_\varepsilon^2\)
\end{equation}
as $\varepsilon\to0$. For any real number $\alpha>1$, we define the parameter set
$$X_\alpha:=Y_\alpha\times\[a_\alpha,\alpha\]\times\[1/\alpha,\alpha\]^k$$
where $a_\alpha:=1/\alpha$ in case $n=4$, $a_\alpha:=-\alpha$ in case $n=5$, and
\begin{equation}\label{Yalpha}
Y_\alpha:=\big\{\tau\in\(T_{\xi_0}M\)^k:\quad\left|\tau_i\right|<\alpha\,\text{ and }\,\left|\tau_i-\tau_j\right|>1/\alpha\quad\forall i\ne j\big\}.
\end{equation}
Here $\left|\cdot\right|$ is the Euclidean norm. As $\alpha\to\infty$, $Y_\alpha$ converges to the set 
\begin{equation}\label{Y}
Y:=\left\{\tau\in\(T_{\xi_0}M\)^k:\,\left|\tau_i-\tau_j\right|\ne0\quad\forall i\ne j\right\}.
\end{equation}
As an easy consequence of \eqref{dist} and the convergence rates of $\mu_{\varepsilon,s}$ and $\delta_\varepsilon$, we obtain that for any $\alpha>1$, there exists $\varepsilon_\alpha\in\(0,\varepsilon_0\)$ such that for any $\varepsilon\in\(0,\varepsilon_\alpha\)$ and $\(\tau,s,t\)\in X_\alpha$, we have $\(z_{\varepsilon,s},\xi_{\varepsilon,\tau},\mu_{\varepsilon,s,t}\)\in\(0,\varepsilon_0\)\times\mathcal{D}_{K_0,\varepsilon_0}$, where $\varepsilon_0$ and $K_0$ are defined by Proposition~\ref{Pr1}. For the sake of simplicity, we denote 
\begin{align*}
&W_{\varepsilon,\tau_i,s,t_i}:=W_{\xi_{\varepsilon,\tau_i},\mu_{\varepsilon,s,t_i}}\,,\quad u_{\varepsilon,\tau,s,t}:=u_{z_{\varepsilon,s},\xi_{\varepsilon,\tau},\mu_{\varepsilon,s,t}}\,,\\
&R_{\varepsilon,\tau,s,t}:=R_{\varepsilon,z_{\varepsilon,s},\xi_{\varepsilon,\tau},\mu_{\varepsilon,s,t}}\,,\quad \text{and}\quad\phi_{\varepsilon,\tau,s,t}:=\phi_\varepsilon\(z_{\varepsilon,s},\xi_{\varepsilon,\tau},\mu_{\varepsilon,s,t}\).
\end{align*}

\smallskip
We state our $C^0$--energy estimates in Proposition~\ref{Pr2} below. We refer to Section~\ref{Sec3} for the proof of this result.

\begin{proposition}\label{Pr2}
We fix $\alpha>0$. As $\varepsilon\to0$, we have
\begin{multline}\label{Energy1}
J_\varepsilon\(u_{\varepsilon,\tau,s,t}+\phi_{\varepsilon,\tau,s,t}\)=kc_0-\frac{e^{-2s/\varepsilon}}{\varepsilon}\Big(c_1\Scal_g\(\xi_0\)s\sum_{i=1}^kt_i^2\\
+c_2\sum_{i=1}^kt_i-c_3\Big)-\frac{e^{-2s/\varepsilon}}{\sqrt\varepsilon}\sum_{i=1}^k\bigg(\frac{c_1}{2}st_i^2D^2\Scal_g\(\xi_0\).\(\tau_i,\tau_i\)\\
+c_4\sum_{j\ne i}\frac{t_it_j}{\left|\tau_i-\tau_j\right|^2}\bigg)+\oo\(\frac{e^{-2s/\varepsilon}}{\sqrt\varepsilon}\)
\end{multline}
in case $n=4$, and 
\begin{multline}\label{Energy2}
J_\varepsilon\(u_{\varepsilon,\tau,s,t}+\phi_{\varepsilon,\tau,s,t}\)=kc_5+c_6\varepsilon^{5/2}-\varepsilon^3\sum_{i=1}^k\Big(c_7\Scal_g\(\xi_0\)t_i^2+c_8t_i^{3/2}\Big)\\
-\varepsilon^{7/2}\Big(c_9s^2+c_8s\sum_{i=1}^kt_i^{3/2}\Big)-\varepsilon^{18/5}\sum_{i=1}^k\bigg(\frac{c_7}{2}t_i^2D^2\Scal_g\(\xi_0\).\(\tau_i,\tau_i\)\\+c_{10}\sum_{j\ne i}\frac{t_i^{3/2}t_j^{3/2}}{\left|\tau_i-\tau_j\right|^3}\bigg)+\oo\(\varepsilon^{18/5}\)
\end{multline}
in case $n=5$, uniformly in $\(\tau,s,t\)\in X_\alpha$, where $c_0,\dotsc,c_{10}$ are positive constants depending only on $\(M,g\)$. 
\end{proposition}

In view of the asymptotic expansions \eqref{Energy1} and \eqref{Energy2}, we introduce the changes of variables
\begin{equation}\label{NewVariables}
\hat{s}=\left\{\begin{aligned}&\varepsilon^{-1/2}\(s-s_0\)&&\text{if }n=4\\&\varepsilon^{-1/20}\(s-s_0\)&&\text{if }n=5\end{aligned}\right.\quad\text{and}\quad \hat{t}=\delta_\varepsilon^{-1}\(t-t_0\),
\end{equation}
where 
$$s_0:=\left\{\begin{aligned}&\frac{c_2}{2c_1\(-\Scal_g\(\xi_0\)\)t_0}&&\text{if }n=4\\&-\frac{kc_8}{2c_9}t_0^{3/2}&&\text{if }n=5\end{aligned}\right.$$
and
$$t_0:=\left\{\begin{aligned}&\frac{2c_3}{kc_2}&&\text{if }n=4\\&\(\frac{3c_8}{4c_7\Scal_g\(\xi_0\)}\)^2&&\text{if }n=5.\end{aligned}\right.$$
Choosing $\alpha$ large enough so that $\(s_0,t_0\)\in\[1/\alpha,\alpha\]^2$ in case $n=4$ and $\(s_0,t_0\)\in\[-\alpha,\alpha\]\times\[1/\alpha,\alpha\]$ in case $n=5$, we can easily see that for any compact subset $A$ of $Y\times\R$, where $Y$ is as in \eqref{Y}, there exists $\varepsilon_A>0$ such that for any $\varepsilon\in\(0,\varepsilon_A\)$, $\(\tau,\hat{s},\hat{t}\)\in A$ implies $\(\tau,s,t\)\in X_\alpha$. Putting together \eqref{Energy1}--\eqref{NewVariables}, we obtain
$$J_\varepsilon\(u_{\varepsilon,\tau,s,t}+\phi_{\varepsilon,\tau,s,t}\)=\left\{\begin{aligned}&kc_0+F_\varepsilon\(\tau,\hat{s},\hat{t}\)\varepsilon^{-1/2}e^{-2s_0/\varepsilon}&&\text{if }n=4\\&kc_5+c_6\varepsilon^{5/2}-\frac{k}{4}c_8t_0^{3/2}\varepsilon^3&&\\&\quad+\frac{k^2c_8^2}{4c_9}t_0^3\varepsilon^{7/2}+F_\varepsilon\(\tau,\hat{s},\hat{t}\)\varepsilon^{18/5}&&\text{if }n=5,\end{aligned}\right.$$
where
\begin{multline}\label{Energy1bis}
F_\varepsilon\(\tau,\hat{s},\hat{t}\)=e^{-2\hat{s}/\sqrt\varepsilon}\bigg(kc_1\(-\Scal_g\(\xi_0\)\)t_0^2\hat{s}+\sum_{i=1}^k\bigg(\frac{c_2}{2t_0}\hat{t}_i^2\\
-\frac{c_1}{2}s_0t_0^2D^2\Scal_g\(\xi_0\).\(\tau_i,\tau_i\)-c_4\sum_{j\ne i}\frac{t_0^2}{\left|\tau_i-\tau_j\right|^2}\bigg)+\oo\(1\)\bigg)
\end{multline}
in case $n=4$, and
\begin{multline}\label{Energy2bis}
F_\varepsilon\(\tau,\hat{s},\hat{t}\)=-c_9\hat{s}^2+\sum_{i=1}^k\bigg(\frac{c_7}{2}\(-\Scal_g\(\xi_0\)\)\hat{t}_i^2\\
-\frac{c_7}{2}t_0^2D^2\Scal_g\(\xi_0\).\(\tau_i,\tau_i\)-c_{10}\sum_{j\ne i}\frac{t_0^{3}}{\left|\tau_i-\tau_j\right|^3}\bigg)+\oo\(1\)
\end{multline}
in case $n=5$, as $\varepsilon\to0$, uniformly in $\(\tau,\hat{s},\hat{t}\)\in A$ for all compact subsets $A$ of $Y\times\R$.

\smallskip
In addition to the above $C^0$--estimates, we need $C^1$--energy estimates in the variables $t_i$. We state these estimates in Proposition~\ref{Pr3} below. We refer to Section~\ref{Sec4} for the proof of this result.

\begin{proposition}\label{Pr3}
Let $A$ be a compact subset of $Y\times\R$. For any $i\in\left\{1,\dotsc,k\right\}$, we have
\begin{equation}\label{Energy3}
\partial_{\hat{t}_i}F_\varepsilon\(\tau,\hat{s},\hat{t}\)=\left\{\begin{aligned}&e^{-2\hat{s}/\sqrt\varepsilon}\(\frac{c_2}{t_0}\hat{t}_i+\oo\(1\)\)&&\text{if }n=4\\&-c_7\Scal_g\(\xi_0\)\hat{t}_i+\oo\(1\)&&\text{if }n=5\end{aligned}\right.
\end{equation}
as $\varepsilon\to0$, uniformly in $\(\tau,\hat{s},\hat{t}\)\in A$. 
\end{proposition}

We are now in position to prove our main result.

\proof[Proof of Theorem~\ref{Th}]
We fix a compact subset $A$ of $Y\times\R$. The choice of $A$ will be precised in the proof. As a consequence of \eqref{Critical}, it suffices to show that for small $\varepsilon>0$, there exists $\(\tau_\varepsilon,\hat{s}_\varepsilon,\hat{t}_\varepsilon\)\in A$ such that $\(z_{\varepsilon,s_\varepsilon},\xi_{\varepsilon,\tau_\varepsilon},\mu_{\varepsilon,s_\varepsilon,t_\varepsilon}\)$ is a critical point of the function $\mathcal{J}_\varepsilon$ defined in Proposition~\ref{Pr1}. As is easily seen, $\(z_{\varepsilon,s_\varepsilon},\xi_{\varepsilon,\tau_\varepsilon},\mu_{\varepsilon,s_\varepsilon,t_\varepsilon}\)$ is a critical point of $\mathcal{J}_\varepsilon$ if and only if $\(\tau_\varepsilon,\hat{s}_\varepsilon,\hat{t}_\varepsilon\)$ is a critical point of $F_\varepsilon$. Here $\hat{s}_\varepsilon$ and $\hat{t}_\varepsilon$ are defined as in \eqref{NewVariables}.

Now we aim to apply Lemma~\ref{CriticalPointLemma} in the appendix to the function $F_\varepsilon$ in a suitable product set. For the sake of clarity, we separate the cases $n=4$ and $n=5$. 

In case $n=4$, we take $A:=\overline{\Omega_1\times\Omega_2}$, where $\Omega_2:=B\(0,r_0\)$ is the open ball in $\R^k$ of center 0 and radius $r_0:=\sqrt{t_0/c_2}$, and $\Omega_1$ is the open subset of $\(T_{\xi_0}M\)^k\times\R$ defined as
$$\Omega_1:=\left\{\(\tau,\hat{s}\)\in Y\times\R:\,G\(\tau\)-1<H\(\hat{s}\)<\inf_YG+1\right\},$$
where
$$G\(\tau\):=\sum_{i=1}^k\bigg(\frac{c_1}{2}s_0t_0^2D^2\Scal_g\(\xi_0\).\(\tau_i,\tau_i\)+c_4\sum_{j\ne i}\frac{t_0^2}{\left|\tau_i-\tau_j\right|^2}\bigg)$$
and
$$H\(\hat{s}\):=kc_1\(-\Scal_g\(\xi_0\)\)t_0^2\hat{s}$$
for all $\(\tau,\hat{s}\)\in Y\times\R$. Since by assumption $D^2\Scal_g\(\xi_0\)$ is definite positive, we obtain that $G>0$ in $Y$ and $A$ is a compact subset of $Y\times\R$. Then Point~(i) in Lemma~\ref{CriticalPointLemma} is an immediate consequence of \eqref{Energy3}. Now we prove Point~(ii). We let $\(\overline\tau,\overline{s}\)\in\Omega_1$ be such that
\begin{equation}\label{MinimumPointCase4}
G\(\overline\tau\)=\inf_YG\quad\text{and}\quad H\(\overline{s}\)=\inf_YG+\frac{1}{2}\,.
\end{equation}
From \eqref{Energy1bis} and \eqref{MinimumPointCase4}, we obtain
\begin{equation}\label{InfCase4}
\inf_{\Omega_2}F_\varepsilon\(\overline\tau,\overline{s},\cdot\)=e^{-2\overline{s}/\sqrt\varepsilon}\(\frac{1}{2}+\oo\(1\)\)
\end{equation}
as $\varepsilon\to0$. By using the fact that $r_0^2<2t_0/c_2$, we also obtain
\begin{equation}\label{SupCase4}
\sup_{\partial\Omega_1\times\Omega_2}F_\varepsilon=\OO\(e^{-2s^*/\sqrt\varepsilon}\)=\oo\(e^{-2\overline{s}/\sqrt\varepsilon}\)
\end{equation}
as $\varepsilon\to0$, where $s^*:=\overline{s}+1/\(2kc_1\(-\Scal_g\(\xi_0\)\)t_0^2\)$ so that $H\(s^*\)=G\(\overline\tau\)+1$. It follows from \eqref{InfCase4} and \eqref{SupCase4} that 
$$\inf_{\Omega_2}F_\varepsilon\(\overline\tau,\overline{s},\cdot\)>\sup_{\partial\Omega_1\times\Omega_2}F_\varepsilon$$
for small $\varepsilon$. Therefore Point~(ii) in Lemma~\ref{CriticalPointLemma} is also satisfied. 

In case $n=5$, we take $A:=\overline{\Omega_1\times\Omega_2}$, , where $\Omega_2:=B\(0,r_0\)$ is the open ball in $\R^k$ of center 0 and radius $r_0:=\sqrt{1/\(-c_7\Scal_g\(\xi_0\)\)}$, and $\Omega_1$ is the open subset of $\(T_{\xi_0}M\)^k\times\R$ defined as
$$\Omega_1:=\left\{\(\tau,\hat{s}\)\in Y\times\R:\,G\(\tau,\hat{s}\)<\inf_YG\(\cdot,0\)+1\right\},$$
where
$$G\(\tau,\hat{s}\):=c_9\hat{s}^2+\sum_{i=1}^k\bigg(\frac{c_7}{2}t_0^2D^2\Scal_g\(\xi_0\).\(\tau_i,\tau_i\)+c_{10}\sum_{j\ne i}\frac{t_0^{3}}{\left|\tau_i-\tau_j\right|^3}\bigg)$$
for all $\(\tau,\hat{s}\)\in Y\times\R$. Similarly to the case $n=4$, we obtain that $G>0$ in $Y\times\R$ and $A$ is a compact subset of $Y\times\R$. Point~(i) in Lemma~\ref{CriticalPointLemma} follows from \eqref{Energy3} together with the assumption that $\Scal_g\(\xi_0\)<0$. To prove Point~(ii), we let $\overline\tau\in Y$ be such that
\begin{equation}\label{MinimumPointCase5}
G\(\overline\tau,0\)=\inf_YG\(\cdot,0\).
\end{equation}
From \eqref{Energy2bis}, \eqref{MinimumPointCase5}, and since $\frac{c_7}{2}\(-\Scal_g\(\xi_0\)\)r_0^2<1$, we obtain
$$\sup_{\partial\Omega_1\times\Omega_2}F_\varepsilon=-\inf_YG\(\cdot,0\)-1+\frac{c_7}{2}\(-\Scal_g\(\xi_0\)\)r_0^2+\oo\(1\)<\inf_{\Omega_2}F_\varepsilon\(0,\overline\tau,\cdot\)$$
for small $\varepsilon$. It follows that Point~(ii) in Lemma~\ref{CriticalPointLemma} is also satisfied.

In both cases $n=4$ and $n=5$, we are now in position to apply Lemma~\ref{CriticalPointLemma} to the function $F_\varepsilon$ in the set $\Omega_1\times\Omega_2$. We obtain that for small $\varepsilon$, there exists a critical point $\(\tau_\varepsilon,\hat{s}_\varepsilon,\hat{t}_\varepsilon\)\in\Omega_1\times\Omega_2$ of $F_\varepsilon$. This ends the proof of 
Theorem~\ref{Th}.
\endproof
 
\section{Proof of the $C^0$--energy estimates}\label{Sec3}

This section is devoted to the proof of Proposition~\ref{Pr2}. We start with proving the following error estimate.

\begin{lemma}\label{ErrorLemma}
We fix $\alpha>0$. We have
\begin{equation}\label{Error0}
R_{\varepsilon,\tau,s,t}=\left\{\begin{aligned}&\OO\(e^{-s/\varepsilon}\)&&\text{if }n=4\\&\OO\(\varepsilon^{9/4}\)&&\text{if }n=5\end{aligned}\right.
\end{equation}
as $\varepsilon\to0$ uniformly in $\(\tau,s,t\)\in X_\alpha$.
\end{lemma}

\proof[Proof of Lemma~\ref{ErrorLemma}]
From the triangular inequality, we obtain
\begin{multline}\label{Error1}
R_{\varepsilon,\tau,s,t}\le\Vol_g\(M\)^{\frac{n+2}{2n}}\left|\varepsilon z_{\varepsilon,s}-z_{\varepsilon,s}^{2^*-1}\right|\\
+\sum_{i=1}^k\left\|\(\Delta_g+\varepsilon\)W_{\varepsilon,\tau_i,s,t_i}-W_{\varepsilon,\tau_i,s,t_i}^{2^*-1}\right\|_{L^{\frac{2n}{n+2}}\(M\)}\\+\Big\|z_{\varepsilon,s}^{2^*-1}+\sum_{i=1}^kW_{\varepsilon,\tau_i,s,t_i}^{2^*-1}-u_{\varepsilon,\tau,s,t}^{2^*-1}\Big\|_{L^{\frac{2n}{n+2}}\(M\)},
\end{multline}
where $\Vol_g\(M\)$ is the volume of the manifold $\(M,g\)$. A straightforward calculation gives
\begin{equation}\label{Error2}
\varepsilon z_{\varepsilon,s}-z_{\varepsilon,s}^{2^*-1}=\left\{\begin{aligned}&\OO\(e^{-s/\varepsilon}\)&&\text{if }n=4\\&\OO\(\varepsilon^{9/4}\)&&\text{if }n=5.\end{aligned}\right.
\end{equation}
For any $i\in\left\{1,\dotsc,k\right\}$, we have (see for instance Robert and V\'etois~\cite{RobVet1})
\begin{align}
&\left\|\(\Delta_g+\varepsilon\)W_{\varepsilon,\tau_i,s,t_i}-W_{\varepsilon,\tau_i,s,t_i}^{2^*-1}\right\|_{L^{\frac{2n}{n+2}}\(M\)}\nonumber\\
&\qquad=\OO\(\mu_{\varepsilon,s}^{\frac{n-2}{2}}\)=\left\{\begin{aligned}&\OO\(e^{-s/\varepsilon}\)&&\text{if }n=4\\&\OO\(\varepsilon^{9/4}\)&&\text{if }n=5.\end{aligned}\right.\label{Error3}
\end{align}
With regard to the last term in the right-hand side of \eqref{Error1}, we have
\begin{multline}\label{Error4}
\Big\|z_{\varepsilon,s}^{2^*-1}+\sum_{i=1}^kW_{\varepsilon,\tau_i,s,t_i}^{2^*-1}-u_{\varepsilon,\tau,s,t}^{2^*-1}\Big\|_{L^{\frac{2n}{n+2}}\(M\)}\\
=\OO\Big(z_{\varepsilon,s}\sum_{i=1}^k\left\|W_{\varepsilon,\tau_i,s,t_i}^{2^*-2}\right\|_{L^{\frac{2n}{n+2}}\(M\)}+z_{\varepsilon,s}^{2^*-2}\sum_{i=1}^k\left\|W_{\varepsilon,\tau_i,s,t_i}\right\|_{L^{\frac{2n}{n+2}}\(M\)}\\
+\sum_{i=1}^k\sum_{j\ne i}\left\|W_{\varepsilon,\tau_i,s,t_i}^{2^*-2}W_{\varepsilon,\tau_j,s,t_j}\right\|_{L^{\frac{2n}{n+2}}\(M\)}\Big).
\end{multline}
Rough estimates give
\begin{align}
\left\|W_{\varepsilon,\tau_i,s,t_i}\right\|_{L^{\frac{2n}{n+2}}\(M\)}&=\OO\(\mu_{\varepsilon,s}^{\frac{n-2}{2}}\),\label{Error5}\\
\left\|W_{\varepsilon,\tau_i,s,t_i}^{2^*-2}\right\|_{L^{\frac{2n}{n+2}}\(M\)}&=\OO\(\mu_{\varepsilon,s}^{\frac{n-2}{2}}\),\label{Error6}
\end{align}
and
\begin{equation}
\left\|W_{\varepsilon,\tau_i,s,t_i}^{2^*-2}W_{\varepsilon,\tau_j,s,t_j}\right\|_{L^{\frac{2n}{n+2}}\(M\)}=\OO\(\mu_{\varepsilon,s}^{n-2}d_g\(\xi_{\varepsilon,\tau_i},\xi_{\varepsilon,\tau_j}\)^{2-n}\)\label{Error7}
\end{equation}
for all $i,j\in\left\{1,\dotsc,k\right\}$, $i\ne j$. The latter estimate can be obtained by splitting the integral into three integrals on the domains\,$M\backslash B\(\xi_0,r_0/2\)$, $B\(\xi_0,r_0/2\)\backslash B\(\xi_{\varepsilon,\tau_i},d_g\(\xi_{\varepsilon,\tau_i},\xi_{\varepsilon,\tau_j}\)/2\)$, and $B\(\xi_{\varepsilon,\tau_i},d_g\(\xi_{\varepsilon,\tau_i},\xi_{\varepsilon,\tau_j}\)/2\)$, and using suitable changes of variable together with the fact that $\mu_{\varepsilon,s}=\oo\(d_g\(\xi_{\varepsilon,\tau_i},\xi_{\varepsilon,\tau_j}\)\)$ as $\varepsilon\to0$. By putting together \eqref{Error4}--\eqref{Error7}, we then obtain
\begin{equation}\label{Error8}
\Big\|z_{\varepsilon,s}^{2^*-1}+\sum_{i=1}^kW_{\varepsilon,\tau_i,s,t_i}^{2^*-1}-u_{\varepsilon,\tau,s,t}^{2^*-1}\Big\|_{L^{\frac{2n}{n+2}}\(M\)}=\left\{\begin{aligned}&\OO\(\varepsilon^{-1}e^{-2s/\varepsilon}\)&&\text{if }n=4\\&\OO\(\varepsilon^3\)&&\text{if }n=5.\end{aligned}\right.
\end{equation}
Finally \eqref{Error0} follows from \eqref{Error1}--\eqref{Error3} and \eqref{Error8}.
\endproof

\proof[Proof of Proposition~\ref{Pr2}]
From \eqref{Energy0} and Lemma~\ref{ErrorLemma}, we obtain
\begin{equation}\label{Energy5}
J_\varepsilon\(u_{\varepsilon,\tau,s,t}+\phi_{\varepsilon,\tau,s,t}\)=J_\varepsilon\(u_{\varepsilon,\tau,s,t}\)+\left\{\begin{aligned}&\OO\(e^{-2s/\varepsilon}\)&&\text{if }n=4\\&\OO\(\varepsilon^{9/2}\)&&\text{if }n=5.\end{aligned}\right.
\end{equation}
Moreover we have
\begin{multline}\label{Energy6}
J_\varepsilon\(u_{\varepsilon,\tau,s,t}\)=J_\varepsilon\(z_{\varepsilon,s}\)+\sum_{i=1}^kJ_\varepsilon\(W_{\varepsilon,\tau_i,s,t_i}\)+\varepsilon z_{\varepsilon,s}\sum_{i=1}^k\int_MW_{\varepsilon,\tau_i,s,t_i}dv_g\\
-z_{\varepsilon,s}\sum_{i=1}^k\int_MW_{\varepsilon,\tau_i,s,t_i}^{2^*-1}dv_g-\frac{1}{2}\sum_{i=1}^k\sum_{j\ne i}\int_MW_{\varepsilon,\tau_i,s,t_i}^{2^*-1}W_{\varepsilon,\tau_j,s,t_j}dv_g\\
+\frac{1}{2}\sum_{i=1}^k\sum_{j\ne i}\int_M\(\Delta_g W_{\varepsilon,\tau_i,s,t_i}+\varepsilon W_{\varepsilon,\tau_i,s,t_i}-W_{\varepsilon,\tau_i,s,t_i}^{2^*-1}\)W_{\varepsilon,\tau_j,s,t_j}dv_g\\
+\frac{1}{2^*}\int_M\bigg(z_{\varepsilon,s}^{2^*}+\sum_{i=1}^kW_{\varepsilon,\tau_i,s,t_i}^{2^*}+2^*z_{\varepsilon,s}\sum_{i=1}^kW_{\varepsilon,\tau_i,s,t_i}^{2^*-1}\\
+2^*\sum_{i=1}^k\sum_{j\ne i}W_{\varepsilon,\tau_i,s,t_i}^{2^*-1}W_{\varepsilon,\tau_j,s,t_j}-u_{\varepsilon,\tau,s,t}^{2^*}\bigg)dv_g\,.
\end{multline}
A straightforward calculation gives
\begin{align}
J_\varepsilon\(z_{\varepsilon,s}\)&=\Vol_g\(M\)\(\frac{\varepsilon z_{\varepsilon,s}^2}{2}-\frac{z_{\varepsilon,s}^{2^*}}{2^*}\)\nonumber\\
&=\left\{\begin{aligned}&\Vol_g\(M\)\frac{e^{-2s/\varepsilon}}{2\varepsilon}+\OO\(\frac{e^{-4s/\varepsilon}}{\varepsilon^4}\)&&\text{if }n=4\\&\Vol_g\(M\)\(\frac{\varepsilon^{5/2}}{5}-\frac{2}{3}s^2\varepsilon^{7/2}\)+\OO\(\varepsilon^4\)&&\text{if }n=5,\end{aligned}\right.\label{Energy7}
\end{align}
where $\Vol_g\(M\)$ is the volume of the manifold $\(M,g\)$. For any $i\in\left\{1,\dotsc,k\right\}$, we have (see for instance Robert and V\'etois~\cite{RobVet1})
\begin{align}
&J_\varepsilon\(W_{\varepsilon,\tau_i,s,t_i}\)=\frac{K_n^{-n}}{n}+\nonumber\\
&\left\{\begin{aligned}&\hspace{-5pt}\frac{K_4^{-4}}{8}\Scal_g\(\xi_{\varepsilon,\tau_i}\)\mu_{\varepsilon,s,t_i}^2\ln\mu_{\varepsilon,s}+\OO\(\mu_{\varepsilon,s}^2+\varepsilon\mu_{\varepsilon,s}^2\left|\ln\mu_{\varepsilon,s}\right|\)\hspace{-4pt}&&\text{if }n=4\\&\hspace{-5pt}-\frac{K_5^{-5}}{10}\Scal_g\(\xi_{\varepsilon,\tau_i}\)\mu_{\varepsilon,s,t_i}^2+\OO\(\mu_{\varepsilon,s}^3+\varepsilon\mu_{\varepsilon,s}^2\)&&\text{if }n=5,\end{aligned}\right.\label{Energy8}
\end{align}
where $K_n$ is the Sobolev constant which was obtained by Rodemich~\cite{Rod}, Aubin~\cite{Aub}, and Talenti~\cite{Tal}, namely
\begin{equation}\label{En8Bis}
\frac{1}{K_n}:=\inf_{u\in D^{1,2}\(\R^n\)\backslash\left\{0\right\}}\frac{\left\|\nabla u\right\|_{L^2\(\R^n\)}}{\left\|u\right\|_{L^{2^*}\(\R^n\)}}=\frac{1}{2}\sqrt{n\(n-2\)}\,\Vol\(\mathbb{S}^n\)^{1/n},
\end{equation}
where $\Vol\(\mathbb{S}^n\)$ is the volume of the standard $n$--dimensional sphere. Moreover since $\xi_0$ is a critical point of $\Scal_g$, a straightforward Taylor expansion gives
\begin{equation}\label{Energy9}
\Scal_g\(\xi_{\varepsilon,\tau_i}\)=\Scal_g\(\xi_0\)+\frac{1}{2}D^2\Scal_g\(\xi_0\)\(\tau_i,\tau_i\)\delta_\varepsilon^2+\OO\(\delta_\varepsilon^3\).
\end{equation}
It follows from \eqref{Energy8} and \eqref{Energy9} that
\begin{multline}
J_\varepsilon\(W_{\varepsilon,\tau_i,s,t_i}\)=\\
\left\{\begin{aligned}&\frac{K_4^{-4}}{4}\bigg(1-\frac{1}{2}\Scal_g\(\xi_0\)st_i^2\frac{e^{-2s/\varepsilon}}{\varepsilon}&&\\
&\quad-\frac{1}{4}D^2\Scal_g\(\xi_0\)\(\tau_i,\tau_i\)st_i^2\frac{e^{-2s/\varepsilon}}{\sqrt\varepsilon}\bigg)+\OO\(\frac{e^{-2s/\varepsilon}}{\varepsilon^{1/4}}\)&&\text{if }n=4\\&\frac{K_5^{-5}}{5}\bigg(1-\frac{1}{2}\Scal_g\(\xi_0\)t_i^2\varepsilon^3&&\\
&\quad-\frac{1}{4}D^2\Scal_g\(\xi_0\)\(\tau_i,\tau_i\)t_i^2\varepsilon^{18/5}\bigg)+\OO\(\varepsilon^{39/10}\)&&\text{if }n=5.\end{aligned}\right.\label{Energy10}
\end{multline}
With regard to the third, fourth, and fifth terms in the right-hand side of \eqref{Energy6}, we obtain
\begin{align}
&\varepsilon z_{\varepsilon,s}\int_MW_{\varepsilon,\tau_i,s,t_i}dv_g=\OO\(\varepsilon z_{\varepsilon,s}\mu_{\varepsilon,s}^{\frac{n-2}{2}}\)=\left\{\begin{aligned}&\OO\(e^{-2s/\varepsilon}\)\hspace{-2pt}&&\text{if }n=4\\&\OO\(\varepsilon^4\)&&\text{if }n=5,\end{aligned}\right.\label{Energy11}\\
&z_{\varepsilon,s}\int_MW_{\varepsilon,\tau_i,s,t_i}^{2^*-1}dv_g=z_{\varepsilon,s}\mu_{\varepsilon,s,t_i}^{\frac{n-2}{2}}\(I_n+\OO\(\mu_{\varepsilon,s}^2|\ln\mu_{\varepsilon,s}|\)\)\nonumber\\
&\quad\qquad\qquad=\left\{\begin{aligned}&\frac{e^{-2s/\varepsilon}}{\varepsilon}t_iI_4+\OO\(\frac{e^{-4s/\varepsilon}}{\varepsilon^2}\)&&\text{if }n=4\\&\(\varepsilon^3+s\varepsilon^{7/2}\)t_i^{3/2}I_5+\OO\(\varepsilon^6|\ln\varepsilon|\)&&\text{if }n=5,\end{aligned}\right.\label{Energy12}
\end{align}
where 
\begin{equation}\label{In}
I_n:=\int_{\R^n}\bigg(\frac{\sqrt{n\(n-2\)}}{1+\left|x\right|^2}\bigg)^{\frac{n+2}{2}}dx\,,
\end{equation}
and
\begin{align}
&\int_MW_{\varepsilon,\tau_i,s,t_i}^{2^*-1}W_{\varepsilon,\tau_j,s,t_j}dv_g\nonumber\\
&\quad=\int_{B\(\xi_{\varepsilon,\tau_i},d_g\(\xi_{\varepsilon,\tau_i},\xi_{\varepsilon,\tau_j}\)/2\)}W_{\varepsilon,\tau_i,s,t_i}^{2^*-1}W_{\varepsilon,\tau_j,s,t_j}dv_g+\OO\(\frac{\mu_{\varepsilon,s}^n}{d_g\(\xi_{\varepsilon,\tau_i},\xi_{\varepsilon,\tau_j}\)^n}\)\nonumber\\
&\quad=\frac{\mu_{\varepsilon,s,t_i}^{\frac{n-2}{2}}\mu_{\varepsilon,s,t_j}^{\frac{n-2}{2}}}{d_g\(\xi_{\varepsilon,\tau_i},\xi_{\varepsilon,\tau_j}\)^{n-2}}\(\[n\(n-2\)\]^{\frac{n-2}{4}}I_n+\oo\(1\)\)\nonumber\\
&\quad=\left\{\begin{aligned}&\frac{e^{-2s/\varepsilon}t_it_j}{\sqrt\varepsilon\left|\tau_i-\tau_j\right|^2}\(2\sqrt2I_4+\oo\(1\)\)&&\text{if }n=4\\&\frac{\varepsilon^{18/5}t_i^{3/2}t_j^{3/2}}{\left|\tau_i-\tau_j\right|^3}\(15^{3/4}I_5+\oo\(1\)\)&&\text{if }n=5\end{aligned}\right.\label{Energy13}
\end{align}
for all $i,j\in\left\{1,\dotsc,k\right\}$, $i\ne j$, where $I_n$ is as in \eqref{In}. To estimate the next term, we observe that 
$$\Delta_g W_{\varepsilon,\tau_i,s,t_i}=W_{\varepsilon,\tau_i,s,t_i}^{2^*-1}+\OO\(W_{\varepsilon,\tau_i,s,t_i}\),$$
which gives
\begin{align}
&\int_M\(\Delta_g W_{\varepsilon,\tau_i,s,t_i}+\varepsilon W_{\varepsilon,\tau_i,s,t_i}-W_{\varepsilon,\tau_i,s,t_i}^{2^*-1}\)W_{\varepsilon,\tau_j,s,t_j}dv_g\nonumber\\
&\qquad=\OO\(\int_MW_{\varepsilon,\tau_i,s,t_i}W_{\varepsilon,\tau_j,s,t_j}dv_g\)\nonumber\\
&\qquad=\left\{\begin{aligned}&\OO\(\mu_{\varepsilon,s}^2\left|\ln\(d_g\(\xi_{\varepsilon,\tau_i},\xi_{\varepsilon,\tau_j}\)\)\right|\)&&\text{if }n=4\\&\OO\big(\mu_{\varepsilon,s}^3d_g\(\xi_{\varepsilon,\tau_i},\xi_{\varepsilon,\tau_j}\)^{-1}\big)&&\text{if }n=5\end{aligned}\right.\nonumber\\
&\qquad=\left\{\begin{aligned}&\OO\(e^{-2s/\varepsilon}\left|\ln\varepsilon\right|\)&&\text{if }n=4\\&\OO\big(\varepsilon^{21/5}\big)&&\text{if }n=5\end{aligned}\right.\label{Energy14}
\end{align}
for all $i,j\in\left\{1,\dotsc,k\right\}$, $i\ne j$. With regard to the last term in the right-hand side of \eqref{Energy6}, we have
\begin{multline}\label{Energy15}
\int_M\hspace{-1pt}\bigg(z_{\varepsilon,s}^{2^*}+\sum_{i=1}^kW_{\varepsilon,\tau_i,s,t_i}^{2^*}+2^*z_{\varepsilon,s}\sum_{i=1}^kW_{\varepsilon,\tau_i,s,t_i}^{2^*-1}\\
+2^*\sum_{i=1}^k\sum_{j\ne i}W_{\varepsilon,\tau_i,s,t_i}^{2^*-1}W_{\varepsilon,\tau_j,s,t_j}-u_{\varepsilon,\tau,s,t}^{2^*}\bigg)dv_g\allowdisplaybreaks\\
=\OO\bigg(z_{\varepsilon,s}^{2^*-1}\sum_{i=1}^k\int_MW_{\varepsilon,\tau_i,s,t_i}dv_g+z_{\varepsilon,s}^2\sum_{i=1}^k\int_MW_{\varepsilon,\tau_i,s,t_i}^{2^*-2}dv_g\\+\sum_{i=1}^k\sum_{j\ne i}\int_MW_{\varepsilon,\tau_i,s,t_i}^{2^*-2}W_{\varepsilon,\tau_j,s,t_j}^2dv_g\bigg).
\end{multline}
Rough estimates give
\begin{align}
\int_MW_{\varepsilon,\tau_i,s,t_i}dv_g&=\OO\(\mu_{\varepsilon,s}^{\frac{n-2}{2}}\),\label{Energy16}\\
\int_MW_{\varepsilon,\tau_i,s,t_i}^{2^*-2}dv_g&=\left\{\begin{aligned}&\OO\(\mu_{\varepsilon,s}^2\left|\ln\mu_{\varepsilon,s}\right|\)\hspace{-2pt}&&\text{if }n=4\\&\OO\(\mu_{\varepsilon,s}^2\)&&\text{if }n=5,\end{aligned}\right.\label{Energy17}
\end{align}
and
\begin{align}
&\int_MW_{\varepsilon,\tau_i,s,t_i}^{2^*-2}W_{\varepsilon,\tau_j,s,t_j}^2dv_g\nonumber\\
&\qquad=\left\{\begin{aligned}&\OO\(\mu_{\varepsilon,s}^4\left|\ln\mu_{\varepsilon,s}\right|d_g\(\xi_{\varepsilon,\tau_i},\xi_{\varepsilon,\tau_j}\)^{-4}\)\hspace{-1pt}&&\text{if }n=4\\&\OO\(\mu_{\varepsilon,s}^4d_g\(\xi_{\varepsilon,\tau_i},\xi_{\varepsilon,\tau_j}\)^{-4}\)&&\text{if }n=5
\end{aligned}\right.\label{Energy18}
\end{align}
for all $i,j\in\left\{1,\dotsc,k\right\}$, $i\ne j$. By combining \eqref{Energy15}--\eqref{Energy18}, we obtain
\begin{align}\label{Energy19}
&\int_M\hspace{-1pt}\bigg(z_{\varepsilon,s}^{2^*}+\sum_{i=1}^kW_{\varepsilon,\tau_i,s,t_i}^{2^*}+2^*z_{\varepsilon,s}\sum_{i=1}^kW_{\varepsilon,\tau_i,s,t_i}^{2^*-1}\nonumber\\
&\quad+2^*\sum_{i=1}^k\sum_{j\ne i}W_{\varepsilon,\tau_i,s,t_i}^{2^*-1}W_{\varepsilon,\tau_j,s,t_j}-u_{\varepsilon,\tau,s,t}^{2^*}\bigg)dv_g\nonumber\\
&\qquad=\left\{\begin{aligned}&\OO\(e^{-2s/\varepsilon}\)&&\text{if }n=4\\&\OO\big(\varepsilon^4\big)&&\text{if }n=5.\end{aligned}\right.
\end{align}
Finally \eqref{Energy1} and \eqref{Energy2} follow from \eqref{Energy5}--\eqref{Energy7}, \eqref{Energy10}--\eqref{Energy14}, and \eqref{Energy19}.
\endproof

\section{Proof of the $C^1$--energy estimates}\label{Sec4}

This section is devoted to the proof of Proposition~\ref{Pr3}.

\proof[Proof of Proposition~\ref{Pr3}]
Throughout this proof, we identify the tangent space $T_\xi M$ with $\mathbb{R}^n$ for all points $\xi$ in a neighborhood of $\xi_0$ by using a smooth, local, orthonormal frame. For any $x\in M$, $\(\tau,s,t\)\in Y\times\R\times\(0,\infty\)^k$, $i\in\left\{1,\dotsc,k\right\}$, and $j\in\left\{1,\dotsc,n\right\}$, we define
$$Z_{\varepsilon,\tau_i,s,t_i,j}\(x\):=\chi\(d_g\(x,\xi_{\varepsilon,\tau_i}\)\)\mu_{\varepsilon,s,t_i}^{\frac{2-n}{2}}V_j\big(\mu_{\varepsilon,s,t_i}^{-1}\exp_{\xi_{\varepsilon,\tau_i}}^{-1}\(x\)\big),$$
where
$$V_0\(y\):=\frac{\left|y\right|^2-1}{\(1+\left|y\right|^2\)^{n/2}}\,\,\text{ and }\,\,V_j\(y\):=\frac{y_j}{\(1+\left|y\right|^2\)^{n/2}}\quad\text{if }j\in\left\{1,\dotsc,n\right\}$$
for all $y\in\R^n$. From Robert and V\'etois~\cite{RobVet2}, we know that the function $\phi_{\varepsilon,\tau,s,t}$ given by Proposition~\ref{Pr1} is such that $\phi_{\varepsilon,\tau,s,t}\in K_{\varepsilon,\tau,s,t}^\perp$ and $DJ_\varepsilon\(u_{\varepsilon,\tau,s,t}+\phi_{\varepsilon,\tau,s,t}\)\in K_{\varepsilon,\tau,s,t}$, where
$$K_{\varepsilon,\tau,s,t}:=\Span\(\left\{1\right\}\cup\left\{Z_{\varepsilon,\tau_i,s,t_i,j}:\,i\in\left\{1,\dotsc,k\right\}\text{ and }j\in\left\{0,\dotsc,n\right\}\right\}\)$$
and
$$K_{\varepsilon,\tau,s,t}^\perp:=\left\{\phi\in H^1\(M\):\quad\left<\phi,\psi\right>_{H^1\(M\)}=0\quad\forall\psi\in K_{\varepsilon,\tau,s,t}\right\}.$$
Let $\lambda_{\varepsilon,\tau,s,t,0}$ and $\lambda_{\varepsilon,\tau,s,t,i,j}$ be real numbers such that
\begin{multline}\label{C1Energy1}
DJ_\varepsilon\(u_{\varepsilon,\tau,s,t}+\phi_{\varepsilon,\tau,s,t}\)=\lambda_{\varepsilon,\tau,s,t,0}\(\frac{d}{d\hat{s}}\[z_{\varepsilon,s}\]\)^{-1}\left<1,\cdot\right>_{H^1\(M\)}\\
+\sum_{i=1}^k\sum_{j=0}^n\lambda_{\varepsilon,\tau,s,t,i,j}\delta_\varepsilon^{-1}\left<Z_{\varepsilon,\tau_i,s,t_i,j},\cdot\right>_{H^1\(M\)}.
\end{multline}
In particular, for any $i_0\in\left\{1,\dotsc,k\right\}$, we obtain
\begin{align}\label{C1Energy2}
&\frac{d}{d\hat{t}_{i_0}}\[J_\varepsilon\(u_{\varepsilon,\tau,s,t}+\phi_{\varepsilon,\tau,s,t}\)\]\nonumber\\
&=\lambda_{\varepsilon,\tau,s,t,0}\(\frac{d}{d\hat{s}}\[z_{\varepsilon,s}\]\)^{-1}\left<1,\frac{d}{d\hat{t}_{i_0}}\[u_{\varepsilon,\tau,s,t}+\phi_{\varepsilon,\tau,s,t}\]\right>_{H^1\(M\)}\nonumber\\
&\,\,\,\,+\sum_{i=1}^k\sum_{j=0}^n\lambda_{\varepsilon,\tau,s,t,i,j}\delta_\varepsilon^{-1}\left<Z_{\varepsilon,\tau_i,s,t_i,j},\frac{d}{d\hat{t}_{i_0}}\[u_{\varepsilon,\tau,s,t}+\phi_{\varepsilon,\tau,s,t}\]\right>_{H^1\(M\)}.
\end{align}
Observe that 
\begin{equation}\label{C1Energy3}
\frac{d}{d\hat{t}_{i_0}}\[u_{\varepsilon,\tau,s,t}\]=\frac{d}{d\hat{t}_{i_0}}\[W_{\varepsilon,\tau_{i_0},s,t_{i_0}}\]=\frac{n^{\frac{n-2}{4}}\(n-2\)^{\frac{n+2}{4}}}{2t_{i_0}}\delta_\varepsilon Z_{\varepsilon,\tau_{i_0},s,t_{i_0},0}\,.
\end{equation}
From now on we fix a compact subset $A$ of $Y\times\R$. All the estimates below will be uniform in $\(\tau,\hat{s},\hat{t}\)\in A$. As $\varepsilon\to0$, rough estimates give
\begin{equation}\label{C1Energy4}
\left<1,Z_{\varepsilon,\tau_{i_1},s,t_{i_1},j_1}\right>_{H^1\(M\)}=\OO\(\mu_{\varepsilon,s}^{\frac{n-2}{2}}\)
\end{equation}
and 
\begin{equation}\label{C1Energy5}
\left<Z_{\varepsilon,\tau_{i_1},s,t_{i_1},j_1},Z_{\varepsilon,\tau_{i_2},s,t_{i_2},j_2}\right>_{H^1\(M\)}=\left\|V_{j_1}\right\|^2_{H^1\(M\)}\delta_{i_1}^{i_2}\delta_{j_1}^{j_2}+\oo\(\delta_\varepsilon\)
\end{equation}
for all $i_1,i_2\in\left\{1,\dotsc,k\right\}$ and $j_1,j_2\in\left\{0,\dotsc,n\right\}$, where $\delta_a^b=0$ if $a\ne b$ and $\delta_a^b=1$ if $a=b$. On the other hand, since $\phi_{\varepsilon,\tau,s,t}\in K_{\varepsilon,\tau,s,t}^\perp$, we obtain
\begin{equation}\label{C1Energy6}
\left<1,\frac{d}{d\hat{t}_{i_0}}\[\phi_{\varepsilon,\tau,s,t}\]\right>_{H^1\(M\)}=0
\end{equation}
and
\begin{align}\label{C1Energy7}
&\left<Z_{\varepsilon,\tau_i,s,t_i,j},\frac{d}{d\hat{t}_{i_0}}\[\phi_{\varepsilon,\tau,s,t}\]\right>_{H^1\(M\)}\nonumber\\
&\qquad=-\left<\frac{d}{d\hat{t}_{i_0}}\[Z_{\varepsilon,\tau_i,s,t_i,j}\],\phi_{\varepsilon,\tau,s,t}\right>_{H^1\(M\)}
\end{align}
for all $i\in\left\{1,\dotsc,k\right\}$ and $j\in\left\{0,\dotsc,n\right\}$. A straightforward computation gives
\begin{equation}\label{C1Energy8}
\left\|\frac{d}{d\hat{t}_{i_0}}\[Z_{\varepsilon,\tau_i,s,t_i,j}\]\right\|_{H^1\(M\)}=\OO\(\delta_\varepsilon\).
\end{equation}
It follows from Cauchy--Schwartz inequality, \eqref{Phi}, \eqref{Error0}, \eqref{C1Energy7}, and \eqref{C1Energy8} that
\begin{equation}\label{C1Energy9}
\left<Z_{\varepsilon,\tau_i,s,t_i,j},\frac{d}{d\hat{t}_{i_0}}\[\phi_{\varepsilon,\tau,s,t}\]\right>_{H^1\(M\)}=\oo\(\delta_\varepsilon\).
\end{equation}
Putting together \eqref{C1Energy2}--\eqref{C1Energy6} and \eqref{C1Energy9}, we obtain
\begin{align}\label{C1Energy10}
&\frac{d}{d\hat{t}_{i_0}}\[J_\varepsilon\(u_{\varepsilon,\tau,s,t}+\phi_{\varepsilon,\tau,s,t}\)\]=\frac{n^{\frac{n-2}{4}}\(n-2\)^{\frac{n+2}{4}}}{2t_{i_0}}\left\|V_0\right\|^2_{H^1\(M\)}\lambda_{\varepsilon,\tau,s,t,i_0,j}\nonumber\\
&\,\,\,\,+\oo\bigg(\(\frac{d}{d\hat{s}}\[z_{\varepsilon,s}\]\)^{-1}\delta_\varepsilon\mu_{\varepsilon,s}^{\frac{n-2}{2}}\left|\lambda_{\varepsilon,\tau,s,t,0}\right|+\delta_\varepsilon\sum_{i=1}^k\sum_{j=0}^n\left|\lambda_{\varepsilon,\tau,s,t,i,j}\right|\bigg).
\end{align}
It remains to estimate the real numbers $\lambda_{\varepsilon,\tau,s,t,0}$ and $\lambda_{\varepsilon,\tau,s,t,i,j}$. We begin with estimating $\lambda_{\varepsilon,\tau,s,t,0}$. From \eqref{C1Energy1} and \eqref{C1Energy4}, we obtain
\begin{multline}\label{C1Energy11}
\lambda_{\varepsilon,\tau,s,t,0}=\Vol_g\(M\)^{-1}DJ_\varepsilon\(u_{\varepsilon,\tau,s,t}+\phi_{\varepsilon,\tau,s,t}\).\frac{d}{d\hat{s}}\[z_{\varepsilon,s}\]\\
+\oo\bigg(\delta_\varepsilon^{-1}\mu_{\varepsilon,s}^{\frac{n-2}{2}}\frac{d}{d\hat{s}}\[z_{\varepsilon,s}\]\sum_{i=1}^k\sum_{j=0}^n\left|\lambda_{\varepsilon,\tau,s,t,i,j}\right|\bigg).
\end{multline}
By observing that 
$$\int_M\phi_{\varepsilon,\tau,s,t}dv_g=\left<1,\phi_{\varepsilon,\tau,s,t}\right>_{H^1\(M\)}=0\,,$$
we obtain
\begin{multline}\label{C1Energy12}
DJ_\varepsilon\(u_{\varepsilon,\tau,s,t}+\phi_{\varepsilon,\tau,s,t}\).1=DJ_\varepsilon\(u_{\varepsilon,\tau,s,t}\).1-\int_M\big[\(u_{\varepsilon,\tau,s,t}+\phi_{\varepsilon,\tau,s,t}\)^{2^*-1}\\
-u_{\varepsilon,\tau,s,t}^{2^*-1}-\(2^*-1\)z_{\varepsilon,s}^{2^*-2}\phi_{\varepsilon,\tau,s,t}\big]dv_g\,.
\end{multline}
Moreover, by using Cauchy--Schwartz and Sobolev inequalities, we obtain
\begin{align}\label{C1Energy13}
&\int_M\big[\(u_{\varepsilon,\tau,s,t}+\phi_{\varepsilon,\tau,s,t}\)^{2^*-1}-u_{\varepsilon,\tau,s,t}^{2^*-1}-\(2^*-1\)z_{\varepsilon,s}^{2^*-2}\phi_{\varepsilon,\tau,s,t}\big]dv_g\nonumber\\
&=\OO\bigg(z_{\varepsilon,s}^{2^*-3}\int_M\phi_{\varepsilon,\tau,s,t}^2dv_g+z_{\varepsilon,s}^{2^*-3}\sum_{i=1}^k\int_M W_{\varepsilon,\tau_i,s,t_i}\left|\phi_{\varepsilon,\tau,s,t}\right|dv_g\nonumber\\
&\quad+\sum_{i=1}^k\int_M W_{\varepsilon,\tau_i,s,t_i}^{2^*-2}\left|\phi_{\varepsilon,\tau,s,t}\right|dv_g+\int_M\left|\phi_{\varepsilon,\tau,s,t}\right|^{2^*-1}dv_g\bigg)\nonumber\\
&=\OO\bigg(z_{\varepsilon,s}^{2^*-3}\left\|\phi_{\varepsilon,\tau,s,t}\right\|^2_{H^1\(M\)}+z_{\varepsilon,s}^{2^*-3}\sum_{i=1}^k\left\|W_{\varepsilon,\tau_i,s,t_i}\right\|_{L^{\frac{2n}{n+2}}\(M\)}\left\|\phi_{\varepsilon,\tau,s,t}\right\|_{H^1\(M\)}\nonumber\\
&\quad+\sum_{i=1}^k\left\|W_{\varepsilon,\tau_i,s,t_i}^{2^*-2}\right\|_{L^{\frac{2n}{n+2}}\(M\)}\left\|\phi_{\varepsilon,\tau,s,t}\right\|_{H^1\(M\)}+\left\|\phi_{\varepsilon,\tau,s,t}\right\|^{2^*-1}_{H^1\(M\)}\bigg).
\end{align}
It follows from \eqref{Phi}, \eqref{Error0}, \eqref{Error5}, \eqref{Error6}, and \eqref{C1Energy13} that
\begin{multline}\label{C1Energy14}
\int_M\big[\(u_{\varepsilon,\tau,s,t}+\phi_{\varepsilon,\tau,s,t}\)^{2^*-1}-u_{\varepsilon,\tau,s,t}^{2^*-1}-\(2^*-1\)z_{\varepsilon,s}^{2^*-2}\phi_{\varepsilon,\tau,s,t}\big]dv_g\\
=\OO\(\mu_{\varepsilon,s}^{n-2}\).
\end{multline}
Putting together \eqref{C1Energy11}, \eqref{C1Energy12}, and \eqref{C1Energy14}, we obtain
\begin{multline}\label{C1Energy15}
\lambda_{\varepsilon,\tau,s,t,0}=\Vol_g\(M\)^{-1}DJ_\varepsilon\(u_{\varepsilon,\tau,s,t}\).\frac{d}{d\hat{s}}\[z_{\varepsilon,s}\]+\OO\(\mu_{\varepsilon,s}^{n-2}\frac{d}{d\hat{s}}\[z_{\varepsilon,s}\]\)\\
+\oo\bigg(\delta_\varepsilon^{-1}\mu_{\varepsilon,s}^{\frac{n-2}{2}}\frac{d}{d\hat{s}}\[z_{\varepsilon,s}\]\sum_{i=1}^k\sum_{j=0}^n\left|\lambda_{\varepsilon,\tau,s,t,i,j}\right|\bigg).
\end{multline}
Now we estimate the real numbers $\lambda_{\varepsilon,\tau,s,t,i,j}$ for all $i\in\left\{1,\dotsc,k\right\}$ and $j\in\left\{0,\dotsc,n\right\}$. From \eqref{C1Energy1}, \eqref{C1Energy4}, and \eqref{C1Energy5}, we obtain
\begin{multline}\label{C1Energy16}
\lambda_{\varepsilon,\tau,s,t,i,j}=\left\|V_j\right\|_{H^1\(M\)}^{-2}DJ_\varepsilon\(u_{\varepsilon,\tau,s,t}+\phi_{\varepsilon,\tau,s,t}\).\(\delta_\varepsilon Z_{\varepsilon,\tau_i,s,t_i,j}\)\\
+\oo\bigg(\(\frac{d}{d\hat{s}}\[z_{\varepsilon,s}\]\)^{-1}\delta_\varepsilon\mu_{\varepsilon,s}^{\frac{n-2}{2}}\left|\lambda_{\varepsilon,\tau,s,t,0}\right|+\delta_\varepsilon\sum_{i'=1}^k\sum_{j'=0}^n\left|\lambda_{\varepsilon,\tau,s,t,i',j'}\right|\bigg).
\end{multline}
By integrating by parts, we obtain
\begin{multline}\label{C1Energy17}
DJ_\varepsilon\(u_{\varepsilon,\tau,s,t}+\phi_{\varepsilon,\tau,s,t}\).Z_{\varepsilon,\tau_i,s,t_i,j}=DJ_\varepsilon\(u_{\varepsilon,\tau,s,t}\).Z_{\varepsilon,\tau_i,s,t_i,j}\\
+\int_M\big[\Delta_gZ_{\varepsilon,\tau_i,s,t_i,j}+\varepsilon Z_{\varepsilon,\tau_i,s,t_i,j}-\(2^*-1\)W_{\varepsilon,\tau_i,s,t_i}^{2^*-2}Z_{\varepsilon,\tau_i,s,t_i,j}\big]\phi_{\varepsilon,\tau,s,t}dv_g\\
-\int_M\big[\(u_{\varepsilon,\tau,s,t}+\phi_{\varepsilon,\tau,s,t}\)^{2^*-1}-u_{\varepsilon,\tau,s,t}^{2^*-1}\\
-\(2^*-1\)W_{\varepsilon,\tau_i,s,t_i}^{2^*-2}\phi_{\varepsilon,\tau,s,t}\big]Z_{\varepsilon,\tau_i,s,t_i,j}dv_g\,.
\end{multline}
By using Cauchy--Schwartz and Sobolev inequalities, we obtain
\begin{multline}\label{C1Energy18}
\int_M\big[\Delta_gZ_{\varepsilon,\tau_i,s,t_i,j}+\varepsilon Z_{\varepsilon,\tau_i,s,t_i,j}-\(2^*-1\)W_{\varepsilon,\tau_i,s,t_i}^{2^*-2}Z_{\varepsilon,\tau_i,s,t_i,j}\big]\phi_{\varepsilon,\tau,s,t}dv_g\\
=\OO\Big(\left\|\Delta_gZ_{\varepsilon,\tau_i,s,t_i,j}+\varepsilon Z_{\varepsilon,\tau_i,s,t_i,j}-\(2^*-1\)W_{\varepsilon,\tau_i,s,t_i}^{2^*-2}Z_{\varepsilon,\tau_i,s,t_i,j}\right\|_{L^{\frac{2n}{n+2}}\(M\)}\\
\times\left\|\phi_{\varepsilon,\tau,s,t}\right\|_{H^1\(M\)}\Big).
\end{multline}
By observing that
$$\Delta_gZ_{\varepsilon,\tau_i,s,t_i,j}+\varepsilon Z_{\varepsilon,\tau_i,s,t_i,j}-\(2^*-1\)W_{\varepsilon,\tau_i,s,t_i}^{2^*-2}Z_{\varepsilon,\tau_i,s,t_i,j}=\OO\(W_{\varepsilon,\tau_i,s,t_i}\),$$
and using \eqref{Phi}, \eqref{Error0}, \eqref{Error5}, and \eqref{C1Energy18}, we obtain
\begin{multline}\label{C1Energy19}
\int_M\big[\Delta_gZ_{\varepsilon,\tau_i,s,t_i,j}+\varepsilon Z_{\varepsilon,\tau_i,s,t_i,j}-\(2^*-1\)W_{\varepsilon,\tau_i,s,t_i}^{2^*-2}Z_{\varepsilon,\tau_i,s,t_i,j}\big]\phi_{\varepsilon,\tau,s,t}dv_g\\
=\OO\(\mu_{\varepsilon,s}^{n-2}\).
\end{multline}
With regard to the last term in the right-hand side of \eqref{C1Energy17}, by observing that $Z_{\varepsilon,\tau_i,s,t_i,j}=\OO\(W_{\varepsilon,\tau_i,s,t_i}\)$ and using Cauchy--Schwartz and Sobolev inequalities, we obtain
\begin{multline}\label{C1Energy20}
\int_M\big[\(u_{\varepsilon,\tau,s,t}+\phi_{\varepsilon,\tau,s,t}\)^{2^*-1}-u_{\varepsilon,\tau,s,t}^{2^*-1}-\(2^*-1\)W_{\varepsilon,\tau_i,s,t_i}^{2^*-2}\phi_{\varepsilon,\tau,s,t}\big]\\
\times Z_{\varepsilon,\tau_i,s,t_i,j}dv_g=\OO\bigg(\int_M\bigg(W_{\varepsilon,\tau_i,s,t_i}^{2^*-3}\left|\phi_{\varepsilon,\tau,s,t}\right|+z_{\varepsilon,s}W_{\varepsilon,\tau_i,s,t_i}^{2^*-3}\\
+\sum_{l\ne i}W_{\varepsilon,\tau_l,s,t_l}W_{\varepsilon,\tau_i,s,t_i}^{2^*-3}+z_{\varepsilon,s}^{2^*-2}+\sum_{l\ne i}W_{\varepsilon,\tau_l,s,t_l}^{2^*-2}+\left|\phi_{\varepsilon,\tau,s,t}\right|^{2^*-2}\bigg)\\
\times W_{\varepsilon,\tau_i,s,t_i}\left|\phi_{\varepsilon,\tau,s,t}\right|dv_g\bigg)=\OO\bigg(\bigg(\left\|W_{\varepsilon,\tau_i,s,t_i}\right\|^{2^*-2}_{H^1\(M\)}\left\|\phi_{\varepsilon,\tau,s,t}\right\|_{H^1\(M\)}\\
+z_{\varepsilon,s}\left\|W_{\varepsilon,\tau_i,s,t_i}^{2^*-2}\right\|_{L^{\frac{2n}{n+2}}\(M\)}+\sum_{l\ne i}\left\|W_{\varepsilon,\tau_l,s,t_l}W_{\varepsilon,\tau_i,s,t_i}^{2^*-2}\right\|_{L^{\frac{2n}{n+2}}\(M\)}\\
+z_{\varepsilon,s}^{2^*-2}\left\|W_{\varepsilon,\tau_i,s,t_i}\right\|_{L^{\frac{2n}{n+2}}\(M\)}+\sum_{l\ne i}\left\|W_{\varepsilon,\tau_l,s,t_l}^{2^*-2}W_{\varepsilon,\tau_i,s,t_i}\right\|_{L^{\frac{2n}{n+2}}\(M\)}\\
\quad+\left\|W_{\varepsilon,\tau_i,s,t_i}\right\|_{H^1\(M\)}\left\|\phi_{\varepsilon,\tau,s,t}\right\|^{2^*-2}_{H^1\(M\)}\bigg)\left\|\phi_{\varepsilon,\tau,s,t}\right\|_{H^1\(M\)}\bigg).
\end{multline}
From \eqref{Phi}, \eqref{Error0}, \eqref{Error5}--\eqref{Error7}, \eqref{C1Energy20}, and since $\left\|W_{\varepsilon,\tau_i,s,t_i}\right\|_{H^1\(M\)}=\OO\(1\)$, we obtain
\begin{multline}\label{C1Energy21}
\int_M\big[\(u_{\varepsilon,\tau,s,t}+\phi_{\varepsilon,\tau,s,t}\)^{2^*-1}-u_{\varepsilon,\tau,s,t}^{2^*-1}-\(2^*-1\)W_{\varepsilon,\tau_i,s,t_i}^{2^*-2}\phi_{\varepsilon,\tau,s,t}\big]\\
\times Z_{\varepsilon,\tau_i,s,t_i,j}dv_g=\OO\(\mu_{\varepsilon,s}^{n-2}\).
\end{multline}
Putting together \eqref{C1Energy16}, \eqref{C1Energy17}, \eqref{C1Energy19}, and \eqref{C1Energy21}, we obtain
\begin{multline}\label{C1Energy22}
\lambda_{\varepsilon,\tau,s,t,i,j}=\left\|V_j\right\|_{H^1\(M\)}^{-2}DJ_\varepsilon\(u_{\varepsilon,\tau,s,t}\).\(\delta_\varepsilon Z_{\varepsilon,\tau_i,s,t_i,j}\)+\OO\(\delta_\varepsilon\mu_{\varepsilon,s}^{n-2}\)\\
+\oo\bigg(\(\frac{d}{d\hat{s}}\[z_{\varepsilon,s}\]\)^{-1}\delta_\varepsilon\mu_{\varepsilon,s}^{\frac{n-2}{2}}\left|\lambda_{\varepsilon,\tau,s,t,0}\right|+\delta_\varepsilon\sum_{i'=1}^k\sum_{j'=0}^n\left|\lambda_{\varepsilon,\tau,s,t,i',j'}\right|\bigg).
\end{multline}
It follows from \eqref{C1Energy3}, \eqref{C1Energy10}, \eqref{C1Energy15}, and \eqref{C1Energy22} that
\begin{multline}\label{C1Energy23}
\frac{d}{d\hat{t}_{i_0}}\[J_\varepsilon\(u_{\varepsilon,\tau,s,t}+\phi_{\varepsilon,\tau,s,t}\)\]=\frac{d}{d\hat{t}_{i_0}}\[J_\varepsilon\(u_{\varepsilon,\tau,s,t}\)\]+\OO\(\delta_\varepsilon\mu_{\varepsilon,s}^{n-2}\)\\
+\oo\bigg(\left|DJ_\varepsilon\(u_{\varepsilon,\tau,s,t}\).1\right|\delta_\varepsilon\mu_{\varepsilon,s}^{\frac{n-2}{2}}\\
+\delta_\varepsilon^2\sum_{i=1}^k\sum_{j=0}^n\left|DJ_\varepsilon\(u_{\varepsilon,\tau,s,t}\).\(\delta_\varepsilon Z_{\varepsilon,\tau_i,s,t_i,j}\)\right|\bigg).
\end{multline}
Similar computations as those performed in Section~\ref{Sec3} give
\begin{align}
&\frac{d}{d\hat{t}_{i_0}}\[J_\varepsilon\(u_{\varepsilon,\tau,s,t}\)\]=\left\{\begin{aligned}&\varepsilon^{-1/2}e^{-2s/\varepsilon}\(\frac{c_2}{t_0}\hat{t}_{i_0}+\oo\(1\)\)&&\text{if }n=4\\&\varepsilon^{18/5}\(-c_7\Scal_g\(\xi_0\)\hat{t}_{i_0}+\oo\(1\)\)&&\text{if }n=5,\end{aligned}\right.\label{C1Energy24}\\
&DJ_\varepsilon\(u_{\varepsilon,\tau,s,t}\).1=\left\{\begin{aligned}&\OO\(e^{-s/\varepsilon}\)&&\text{if }n=4\\&\OO\(\varepsilon^{23/10}\)&&\text{if }n=5,\end{aligned}\right.\label{C1Energy25}
\end{align}
and
\begin{equation}\label{C1Energy26}
DJ_\varepsilon\(u_{\varepsilon,\tau,s,t}\).Z_{\varepsilon,\tau_i,s,t_i,j}=\left\{\begin{aligned}&\OO\(\varepsilon^{-1}e^{-2s/\varepsilon}\)&&\text{if }n=4\\&\OO\(\varepsilon^3\)&&\text{if }n=5\end{aligned}\right.
\end{equation}
for all $i\in\left\{1,\dotsc,k\right\}$ and $j\in\left\{0,\dotsc,n\right\}$. From \eqref{C1Energy23}--\eqref{C1Energy26}, we obtain
$$\frac{d}{d\hat{t}_{i_0}}\[J_\varepsilon\(u_{\varepsilon,\tau,s,t}+\phi_{\varepsilon,\tau,s,t}\)\]=\left\{\begin{aligned}&\varepsilon^{-1/2}e^{-2s/\varepsilon}\(\frac{c_2}{t_0}\hat{t}_{i_0}+\oo\(1\)\)&&\text{if }n=4\\&\varepsilon^{18/5}\(-c_7\Scal_g\(\xi_0\)\hat{t}_{i_0}+\oo\(1\)\)&&\text{if }n=5.\end{aligned}\right.$$
This ends the proof of Proposition~\ref{Pr3}.
\endproof

\appendix

\section{A critical point result for product sets}\label{Appendix}

In this appendix, we prove a critical point result which was used in Section~\ref{Sec2}. This result relies on a deformation argument using a negative gradient-type flow. A similar argument was used by Chen, Wei, and Yan~\cite{ChenWeiYan} in the case of a function of two real variables. 

\medskip
The Lyapunov--Schmidt method crucially depends on the existence of critical points for families of functions $\(F_\varepsilon\)_{\varepsilon>0}$ which converge to a function $F_0$. In case the limit function $F_0$ has a saddle point $x_0$, if the functions $F_\varepsilon$ converge only in $C^0$ to $F_0$, then it is not true in general that there exist critical points of the functions $F_\varepsilon$ which converge to $x_0$, even when assuming that $x_0$ is a non-degenerate critical point of $F_0$. From degree theory, we know that this property holds true if we replace $C^0$--convergence by $C^1$--convergence and we assume that the critical point $x_0$ is non-degenerate. The objective of the result below is to obtain this property under weaker conditions which only involve derivatives in some directions.


\begin{lemma}\label{CriticalPointLemma}
Let $n_1,n_2\ge1$ be two integers, $\Omega_1$ be a bounded and open subset of $\R^{n_1}$, $\Omega_2$ be a bounded, open, and smooth subset of $\R^{n_2}$, and $\Omega:=\Omega_1\times\Omega_2$.~Let $F$ be a $C^2$--function in a neighborhood of $\overline\Omega$ such that
\begin{enumerate}
\renewcommand{\labelenumi}{(\roman{enumi})}
\item The outward normal derivative of $F$ on $\Omega_1\times\partial\Omega_2$ is positive.
\item There exists $\overline{x}\in\Omega_1$ such that $\displaystyle\inf_{\Omega_2}F\(\overline{x},\cdot\)>\sup_{\partial\Omega_1\times\Omega_2}F$.
\end{enumerate}
Then $F$ has a critical point in (the interior of) $\Omega$. 
\end{lemma}

\proof[Proof of Lemma~\ref{CriticalPointLemma}]
We assume by contradiction that the function $F$ does not have any critical point in $\Omega$.

We start our proof by constructing a negative gradient-type flow for the function $F$. From Point~(ii) and the continuity of $F$ on $\overline\Omega$, we obtain that there exists an open set $U$ such that $\overline U\subset\Omega_1$ and
\begin{equation}\label{CPLem1}
\inf_{\Omega_2}F\(\overline{x},\cdot\)>\sup_{\(\Omega_1\backslash U\)\times\Omega_2}F.
\end{equation}
We let $V$ and $W$ be two open sets such that $\overline U\subset V$, $\overline V\subset W$, and $\overline W\subset\Omega_1$. We let $\chi$ be a smooth cutoff function in $\R^{n_1}$ such that $\chi\equiv1$ in $V$, $\chi\equiv0$ in $\R^{n_1}\backslash W$, and $0\le\chi\le1$ in $W\backslash V$. For $i\in\left\{1,2\right\}$, we let $p_i:\R^{n_1+n_2}\to\R^{n_i}$ be the canonical projection, namely
$$p_i\(x_1,x_2\):=x_i\quad\forall x_1\in\R^{n_1},\,x_2\in\R^{n_2}\,.$$
By assumption, we have that there exists an open subset $D$ of $\R^{n_1+n_2}$ such that $F\in C^2\(D\)$ and $\overline\Omega\subset D$. From basic theory of ODEs, we then obtain the existence of a lower semi-continuous mapping $T:D\mapsto\(0,\infty\]$ and a $C^2$--mapping $\Phi:D_T\mapsto\R^{n_1+n_2}$, where $D_T:=\left\{\(t,x\):\,x\in D\text{ and }t\in\[0,T\(x\)\)\right\}$, such that for any $x\in D$, we have
$$\left\{\begin{aligned}&\frac{\partial\Phi}{\partial t}\(t,x\)=-\chi\(p_1\(\Phi\(t,x\)\)\)\nabla F\(\Phi\(t,x\)\)&&\forall t\in\[0,T\(x\)\)\\&\phi\(0,x\)=x\end{aligned}\right.$$
and either $T\(x\)=\infty$ or $\Phi\(t,x\)\not\in\overline{\Omega}$ when $t$ approaches $T\(x\)$.

We prove that $\Phi\(t,x\)\in\overline\Omega$ for all $x\in\overline\Omega$ and $t\in\[0,T\(x\)\)$, which implies in particular $T\(x\)=\infty$. We assume by contradiction that the curve $t\mapsto\Phi\(t,x\)$ leaves the set $\overline\Omega$, namely that there exist $t_-,t_+\in\[0,T\(x\)\)$ such that $t_-<t_+$, $\Phi\(t_-,x\)\in\partial\Omega$ and $\Phi\(t,x\)\not\in\overline{\Omega}$ for all $t\in\(t_-,t_+\)$. Since $\Phi\(t,x\)$ is not constant in $t$, we infer from the uniqueness of the flow that $\frac{\partial\Phi}{\partial t}\(t_-,x\)\ne0$. It follows that $\chi\(p_1\(\Phi\(t_-,x\)\)\)\ne0$, which gives $\Phi\(t_-,x\)\in\Omega_1\times\partial\Omega_2$. From Point~(i), we then obtain
\begin{equation}\label{CPLem2}
\frac{d}{dt}\<\Phi\(t_-,x\),\nu\>=-\chi\(p_1\(\Phi\(t_-,x\)\)\)\<\nabla F\(\Phi\(t_-,x\)\),\nu\><0,
\end{equation}
where $\<\cdot,\cdot\>$ is the Euclidean scalar product and $\nu$ is the outward normal vector to $\Omega_1\times\partial\Omega_2$ at the point $\Phi\(t_-,x\)$. This contradicts the fact that $\Phi\(t,x\)\not\in\overline{\Omega}$ for all $t\in\(t_-,t_+\)$. Therefore we have proven that $T\(x\)=\infty$ and $\Phi\(t,x\)\in\overline\Omega$ for all $x\in\overline\Omega$ and $t\in\[0,\infty\)$.

Now we define
$$c:=\inf_{h\in\Gamma}\sup_{x\in\Omega}F\(h\(x\)\),$$
where
$$\Gamma:=\left\{h\in C^0\(\overline\Omega,\overline\Omega\):\,h\(x\)=x\quad\forall x\in\partial\Omega_1\times\Omega_2\right\}.$$
Our aim is to construct a mapping $h_0\in\Gamma$ such that 
\begin{equation}\label{CPLem3}
\sup_{x\in\Omega}F\(h_0\(x\)\)<c
\end{equation}
so to obtain a contradiction. 

Since $\overline U\subset V$ and $\Phi\in C^0\(\[0,\infty\)\times\Omega,\Omega\)$, we obtain that there exists a real number $t_0>0$ such that $\Phi\(t,U\times\Omega_2\)\subset V\times\Omega_2$ for all $t\in\[0,t_0\]$. Since $F\in C^1\(\overline\Omega\)$, $\overline{V}\subset\Omega_1$, $\nabla F\ne0$ on $\Omega_1\times\partial\Omega_2$ according to Point~(i), and we have assumed at the beginning of the proof that $\nabla F\ne0$ in $\Omega$, we obtain the existence of a real number $\delta_0>0$ such that $\left|\nabla F\right|\ge\delta_0$ in $V\times\Omega_2$. From the definition of $c$, we obtain that there exists $h\in\Gamma$ such that
\begin{equation}\label{CPLem4}
\sup_{x\in\Omega}F\(h\(x\)\)\le c+\frac{t_0\delta_0^2}{2}\,.
\end{equation}
Now we define $h_0:=\Phi\(t_0,h\)$, and we prove \eqref{CPLem3}. We separate the cases $h\(x\)\in U\times\Omega_2$ and $h\(x\)\in\(\Omega_1\backslash U\)\times\Omega_2$. In case $h\(x\)\in U\times\Omega_2$, since $\Phi\(t,U\times\Omega_2\)\subset V\times\Omega_2$ for all $t\in\[0,t_0\]$, $\chi\equiv1$ in $V$, and $\left|\nabla F\right|\ge\delta_0$ in $V\times\Omega_2$, we obtain
\begin{equation}\label{CPLem5}
F\(h\(x\)\)-F\(h_0\(x\)\)=\int_0^{t_0}\left|\nabla F\(\Phi\(t,h\(x\)\)\)\right|^2dt\ge t_0\delta_0^2\,.
\end{equation}
It follows from \eqref{CPLem4} and \eqref{CPLem5} that
\begin{equation}\label{CPLem6}
\sup_{x\in h^{-1}\(U\times\Omega_2\)}F\(h_0\(x\)\)\le c-\frac{t_0\delta_0^2}{2}\,.
\end{equation}
On the other hand, since the function $t\mapsto F\(\Phi\(t,h\(x\)\)\)$ is nonincreasing for all $x\in h^{-1}\(\(\Omega_1\backslash U\)\times\Omega_2\)$, it follows from \eqref{CPLem1} that
\begin{equation}\label{CPLem7}
\sup_{x\in h^{-1}\(\(\Omega_1\backslash U\)\times\Omega_2\)}F\(h_0\(x\)\)<\inf_{\Omega_2}F\(\overline{x},\cdot\).
\end{equation}
It remains to prove
\begin{equation}\label{CPLem8}
\inf_{\Omega_2}F\(\overline{x},\cdot\)\le c\,.
\end{equation}
We fix a point $\overline{y}\in\Omega_2$. For any mapping $h\in\Gamma$, we define $\overline{h}:=p_1\(h\(\cdot,\overline{y}\)\)$. We infer from the properties of $h$ that $\overline{h}\in C^0\(\overline{\Omega_1},\overline{\Omega_1}\)$ and $\overline{h}\(x\)=x$ for all points $x\in\partial\Omega_1$. We then obtain from degree theory that $\overline{h}\(\overline{\Omega_1}\)=\overline{\Omega_1}$ (see Poincar\'e--Bohl theorem in~\cite{OrtRhe}). In particular, we obtain that there exists a point $x_0\in\Omega_1$ such that $\overline{h}\(x_0\)=\overline{x}$. From the definition of $\overline{h}$, it follows that there exists a point $y_0\in\Omega_2$ such that $h\(x_0,\overline{y}\)=\(\overline{x},y_0\)$. We then obtain
\begin{equation}\label{CPLem9}
\inf_{\Omega_2}F\(\overline{x},\cdot\)\le F\(\overline{x},y_0\)=F\(h\(x_0,\overline{y}\)\)\le\sup_{x\in\Omega}F\(h\(x\)\).
\end{equation}
Since \eqref{CPLem9} holds true for all mappings $h\in\Gamma$, we obtain \eqref{CPLem8}. 

Finally \eqref{CPLem3} follows from \eqref{CPLem6}, \eqref{CPLem7}, and \eqref{CPLem8}. This ends the proof of Lemma~\ref{CriticalPointLemma}.
\endproof

\section*{Acknowledgments} This work has been initiated ans partially carried out during the visit of the first author at McGill University. The first author gratefully acknowledges the hospitality and the financial support of this institution during his stay.

\end{document}